\documentclass[reqno,12pt,a4]{amsart}
\usepackage{amssymb}
\usepackage{amsmath}

\def\Ext{\mathop{\rm Ext}\nolimits}

\def\supp{\mathop{\rm supp}\nolimits}

\def\Op{\mathop{\rm Op}\nolimits}

\newtheorem{thm}{Theorem}
\newtheorem{lem}[thm]{Lemma}

\newtheorem{dfn}[thm]{Definition}

\def\K{{\mathcal K}}

\date{}
\author{V. Manuilov}\thanks{Partially supported by
RFFI grant No. 02-01-00574 and by HIII\hspace{-2.3ex}%
\rule{1.9ex}{0.07ex}\,-619.2003.01.}
\title[Equivalence of two approaches to index theory]{Translation
invariant asymptotic homomorphisms: equivalence of two approaches
in the index theory}

\begin{document}

\begin{abstract}
The algebra $\Psi(M)$ of order zero pseudodifferential
operators on a compact manifold $M$ defines a well-known $C^*$-extension
of the algebra $C(S^*M)$ of continuous functions on the cospherical bundle
$S^*M\subset T^*M$ by the algebra $\K$ of compact operators.
In his proof of the index theorem, Higson defined and used an
asymptotic homomorphism $T$ from $C_0(T^*M)$ to $\K$, which plays the role
of a deformation for the commutative algebra $C_0(T^*M)$.
Similar constructions exist also for operators and
symbols with coefficients in a $C^*$-algebra. We show that the image of
the above extension under the Connes--Higson construction is $T$
and that this extension can be reconstructed out of $T$. This explains, why
the classical approach to the index theory coincides with the one based on
asymptotic homomorphisms.

\end{abstract}

\maketitle

\section{Two ways to define index}

The standard way to define the index of a pseudodifferential elliptic
operator on a compact manifold $M$ comes from the short exact sequence
of $C^*$-algebras
 \begin{equation}\label{extens}
0\to \K\to \Psi(M)\to C(S^*M)\to 0,
 \end{equation}
where $\K$ is the algebra of compact operators on $L^2(M)$,
$\Psi(M)$ denotes the norm closure of the algebra of order zero
pseudodifferential operators in the algebra of bounded operators
on $L^2(M)$ and $S^*M$ denotes the cospherical bundle, $S^*M=\{(x,\xi)\in
T^*M:|\xi|=1\}$, in the cotangent bundle $T^*M$. If one deals with
operators having coefficients in a $C^*$-algebra $A$ then one has to
tensor the short exact sequence (\ref{extens}) by $A$:
 \begin{equation}\label{extensi}
0\to \K\otimes A\to \Psi_A(M)\to C(S^*M;A)\to 0,
 \end{equation}
where $C(X;A)$ denotes the $C^*$-algebra of continuous functions on $X$
taking values in $A$.
The (main) symbol of a pseudodifferential elliptic operator of order zero
is an invertible element in $C(S^*M;A)$ and the $K$-theory boundary
map $K_1(C(S^*M;A))\to K_0(\K\otimes A)$ maps the symbol to a class in
$K_0(A)$, which is called the index of the operator.

Another approach, suggested by Higson in \cite{Hig}, is based on the
notion of an asymptotic homomorphism \cite{CH}.
Here one starts with a symbol
$\sigma$ of a pseudodifferential operator of order one and constructs
a {\it symbol class} $[a_\sigma]\in K_0(C_0(T^*M))$ (see details in
\cite{Hig}).
Then one constructs an asymptotic homomorphism from $C_0(T^*M)$ to $\K$ as
follows. In the local coordinates $(x,\xi)$ in $U\times\mathbb R^n\subset
TM$ take a smooth function $a(x,\xi)$ with a compact support,
$a\in C_c^\infty(U\otimes\mathbb R^n)$. Then define a
continuous family of operators $T_{a,t}$, $t\in\mathbb R_+=(0,\infty)$,
on $L^2(U)$ by
 \begin{equation}\label{asymT}
T_{a,t}f(x)=\int a(x,t^{-1}\xi)e^{ix\xi}\widehat{f}(\xi)d\xi,
 \end{equation}
where $\widehat{f}$ is the
Fourier transform for $f$. Fix an atlas $\{U_k\}$ of charts on a compact
manifold $M$ and let $\{\varphi_k\}$ be a smooth partition of unity
subordinate to that atlas. Take also smooth functions $\psi_k$ on $M$ such
that $\supp\psi_k\subset U_k$ and $\psi_k\varphi_k=\varphi_k$ for all $k$.
For $f\in L^2(M)$ put
 \begin{equation}\label{asympT}
T_t(a)f=\sum_k T_{\psi_ka,t}(\varphi_kf).
 \end{equation}
It is shown in \cite{Hig}, Lemma 8.4, that this family of operators
defines an asymptotic homomorphism $T=(T_t)_{t\in\mathbb R_+}$ from
$C_c^\infty(T^*M)$ to $\K$ (it is also shown in \cite{Hig}, Lemma 8.7,
that if we take another atlas or other
functions $\varphi_k$ and $\psi_k$ then the resulting asymptotic
homomorphism is asymptotically equal to this one). Therefore this
asymptotic homomorphism defines a $*$-homomorphism $\overline{T}$ from
$C_c^\infty(T^*M)$ to
the asymptotic $C^*$-algebra $C_b(\mathbb R_+;\K)/C_0(\mathbb R_+;\K)$,
where $C_b(\mathbb R_+;\K)$ denotes the algebra of bounded continuous
$\K$-valued functions on $\mathbb R_+$ and due to automatic continuity of
$C^*$-algebra
$*$-homomorphisms one can extend $\overline{T}$ to a $*$-homomorphism
$\widehat{T}:C_0(T^*M)\to C_b(\mathbb R_+;\K)/C_0(\mathbb R_+;\K)$.
Applying the Bartle--Graves selection theorem \cite{Loring}, we obtain an
asymptotic homomorphism
$\widetilde{T}=(\widetilde{T}_t)_{t\in\mathbb R_+}:C_0(T^*M)\to\K$, which
is uniformly continuous and
asymptotically equal to $T$ on smooth functions
(i.e. $\lim_{t\to\infty}T_t(a)-\widetilde{T}_t(a)=0$ for any $a\in
C_c^\infty(T^*M)$). Finally the index of the operator with a
symbol $\sigma$
is defined as the class of the image of $[a_\sigma]$
under the map $K_0(C_0(T^*M))\to K_0(\K)$ induced by $\widetilde{T}$.
Once more,
one can tensor everything by $A$ and construct an asymptotic homomorphism
 $$
T=(T_t)_{t\in\mathbb R_+}:C_c^\infty(T^*M;A)\to\K\otimes A
 $$
and then change it by a uniformly (with respect to $t$) continuous
asymptotic homomorphism extended to $C_0(T^*M;A)$,
 \begin{equation}\label{asympTT}
\widetilde{T}=(\widetilde{T}_t)_{t\in\mathbb R_+}:C_0(T^*M;A)\to
\K\otimes A
 \end{equation}
(we keep the same notation $T$ and $\widetilde{T}$ for the case
of $A$-valued symbols).
Remark that the asymptotic homomorphism $T$ is translation invariant, i.e.
$T_{ts}(a)=T_t(a_s)$, where $a_s(x,\xi)= a(x,s^{-1}\xi)$, $a\in
C_c^\infty(T^*M;A)$. The asymptotic homomorphism $\widetilde{T}$ enjoys
the property of asymptotic translation invariance, i.e.
$\lim_{t\to\infty}\widetilde{T}_{ts}(a)-\widetilde{T}_t(a_s)=0$ for any
$a\in C_0(T^*M;A)$.

The purpose of the present paper is to explain, why these two approaches
to define the index of elliptic pseudodifferential operators are
equivalent.

\section{The Connes--Higson construction and its inverse}

In the pioneering paper on asymptotic homomorphisms, \cite{CH},
a construction was given, which
transforms $C^*$-algebra extensions into asymptotic homomorphisms.
Given a $C^*$-extension $0\to B\to E\to D\to 0$, one obtains an asymptotic
homomorphism from the suspension $SD=C_0((0,1);D)$ into $B$. One of
the main results of \cite{MT6} was that the asymptotic homomorphisms
obtained via this Connes--Higson construction possess an additional
important property --- translation invariance. In order to make its
description easier we identify suspension $SA$ with $C_0(\mathbb
R_+;D)$ instead of using $(0,1)$. There is a natural action
$\tau$ of $\mathbb R_+$ on itself by multiplication, $\tau_s(x)=xs$,
$s,x\in\mathbb R_+$, which extends to an action on $SD$ by
$\tau_s(f)(x)=f(sx)$, where $f\in SD=C_0(\mathbb R_+;D)$
(in \cite{MT6} the additive structure on
$\mathbb R$ was used instead of the multiplicative structure on $\mathbb
R_+$).

\begin{dfn}
{\rm
An asymptotic
homomorphism $\varphi=(\varphi_t)_{t\in\mathbb R_+}:SD\to B$ is
{\it translation invariant} if
$\varphi_t(\tau_s(f))=\varphi_{ts}(f)$ for any $f\in
SD$ and for any $t,s\in\mathbb R_+$ and if $\lim_{t\to 0}\varphi_t(f)=0$
for any $f\in SD$. It is {\it asymptotically translation invariant} if
$\lim_{t\to\infty}\varphi_t(\tau_s(f))-\varphi_{ts}(f)=0$ for any $f\in
SA$ and for any $t,s\in\mathbb R_+$ and if $\lim_{t\to 0}\varphi_t(f)=0$
for any $f\in SD$.
 Two (asymptotically) translation invariant asymptotic homomorphisms
$\varphi^{(0)},\varphi^{(1)}:SD\to B$ are homotopic if there is an
(asymptotically) translation invariant asymptotic homomorphism $\Phi:SD\to
C[0,1]\otimes B$, whose restrictions onto the endpoints of $[0,1]$
coincide with $\varphi^{(0)}$ and $\varphi^{(1)}$ respectively.
}
\end{dfn}

Note that, by passing to spherical coordinates in the fibers, the
suspension $SC(S^*M;A)$ can be identified with the algebra
 $C_{00}(T^*M;A)$ of continuous functions on $T^*M$ vanishing both at
 infinity and at the zero section and the asymptotic homomorphism
$\widetilde{T}$ (\ref{asympTT}) can be restricted onto $C_{00}(T^*M;A)$.

\begin{lem}
The asymptotic homomorphism $\widetilde{T}:C_{00}(T^*M;A)\to
\K\otimes A$ is
asymptotically translation invariant.
\end{lem}
\begin{proof}
The family of maps (\ref{asympTT}) is obviously asymptotically invariant
under the action of $\mathbb R_+$ and one easily checks that $\lim_{t\to
0}\widetilde{T}_t(a)=0$ for any $a\in C_{00}(T^*M;A)$.
\end{proof}

From now on we assume that $A$ and $D$ are separable and $B$ is stable and
$\sigma$-unital.
Let $\Ext_h(D,B)$ denote the semigroup of homotopy classes of
$C^*$-extensions of $D$ by $B$ and let
$[[SD,B]]_{a,\tau}$ denote the semigroup of asymptotically translation
invariant asymptotic homomorphisms from $SD$ to $B$.
Note that there is a forgetful map
 \begin{equation}\label{forge}
[[SD,B]]_{a,\tau}\to[[SD,B]]
 \end{equation}
to the group of homotopy
classes of all asymptotic homomorphisms from $SD$ to $B$, which is the
$E$-theory group $E(SD,B)$. The Connes--Higson construction \cite{CH}
defines a map
 $$
CH:\Ext_h(D,B)\to [[SD,B]].
 $$
It was shown in
\cite{MT6} that this map factorizes through the map (\ref{forge})
and the modified Connes--Higson construction
 \begin{equation}\label{CH}
\widetilde{CH}:\Ext_h(D,B)\to [[SD,B]]_{a,\tau}.
 \end{equation}

The main result of \cite{MT6} is that the map (\ref{CH}) is an
isomorphism.
This was proved by constructing an inverse map
 $$
I:[[SD,B]]_{a,\tau}\to \Ext_h(D,B).
 $$
The map $I$ is constructed as follows (see details in \cite{MT6}).
Let $\varphi=(\varphi_t)_{t\in\mathbb R_+}:SD\to B$ be an asymptotically
translation invariant asymptotic homomorphism. Then, by the Bartle--Graves
continuous selection theorem \cite{Loring}, there is an asymptotically
translation invariant asymptotic homomorphism $\widetilde{\varphi}$, which
is asymptotically equal to $\varphi$ and such that the family of maps
$\widetilde{\varphi}_t:SD\to B$ is uniformly continuous.

Let $\gamma_0\in C_0(\mathbb R_+)$ be a (smooth)
function with support in $[1/2,2]$ such that
$\sum_{i\in\mathbb Z}\gamma_i^2=1$, where $\gamma_i=\tau_{2^i}(\gamma_0)$.
Note that $\gamma_i\gamma_j=0$ when $|i-j|\geq 2$.
Let $e_{ij}$ denote the standard elementary operators on the standard
Hilbert $C^*$-module $H_B=l^2(\mathbb Z)\otimes B$.
We identify the algebra of compact (resp.
adjointable) operators on $H_B$ with the $C^*$-algebra $B\otimes\K$ (resp.
the multiplier $C^*$-algebra
$M(B\otimes\K)$) and let
 $$
q:M(B\otimes\K)\to
Q(B\otimes\K)=M(B\otimes\K)/B\otimes\K
 $$
be the quotient $*$-homomorphism.

Put, for $a\in D$,
 $$
I_0(\varphi)(a)=\sum_{i,j\in\mathbb Z}\widetilde{\varphi}_{2^i}
(\tau_{2^{-i}}(\gamma_i\gamma_j)\otimes a)\otimes e_{ij}\in M(B\otimes\K)
 $$
and $I(\varphi)(a)=q(I_0(\varphi)(a))$.
The map $I:D\to Q(B\otimes\K)$ is a $*$-homomorphism, so it defines an
extension of $D$ by $B\otimes\K$, being its Busby invariant.

\section{Main result}

Denote by $[\Psi_A(M)]\in\Ext_h(C(S^*M;A),\K\otimes A)$ the homotopy
class of the extension (\ref{extensi}).

\begin{thm}
The image of $[\Psi_A(M)]$ under the Connes--Higson construction coincides
with the homotopy class of the asymptotic homomorphism
$\widetilde{T}$ if $A$ is separable.
\end{thm}
\begin{proof}
Due to \cite{MT6} we do not need to prove that $CH([\Psi_A(M)])$ is
homotopic to $\widetilde{T}$.
It is sufficient to prove instead that $I(\widetilde{T})$ is
homotopic to the Busby invariant of the extension (\ref{extensi}), which
is easier.

In order to construct the Busby invariant for the extension (\ref{extensi})
one can use the same atlas of charts and the same functions $\varphi_k$ and
$\psi_k$ as in the construction of the asymptotic homomorphism $T_t$
(\ref{asympT}).
Let $\theta$ be a smooth cutting function on $[0,\infty)$, which equals 1
outside a compact set and vanishes at the origin and let
$U\subset M$ be a subset diffeomorphic to a domain in
a Euclidean space. In the local coordinates $(x,\xi)$ in
$U\times\mathbb R^n\subset T^*M$ take a smooth function
$a(x,\xi)$ with a compact support with respect to the first coordinate
and order zero homogeneous with respect to the second coordinate.
Let $f$ be an element of the
Hilbert $C^*$-module $L^2(U)\otimes A$ over $A$. Define an operator
$\Op(a)$ on this Hilbert $C^*$-module by
 $$
\Op(a)f(x)=\int a(x,\xi)\theta(|\xi|)e^{ix\xi}\widehat{f}(\xi)\,d\xi,
 $$
where $\widehat{f}$ is the Fourier transform for $f$.
Then, for a main symbol $a(x,\xi)\in
C^\infty(S^*M;A)$ defined on the whole $M$,
one can construct an operator $\Op(a)$ on the Hilbert $C^*$-module
$L^2(M)\otimes A$ by
 $$
\Op(a)(f)=\sum_k\Op(\psi_ka)(\varphi_kf),\qquad
\Op(a)\in M(\K\otimes A).
 $$
The map $q\circ\Op:C^\infty(S^*M;A)\to Q(\K\otimes A)$
is a $*$-homomorphism (cf.
\cite{Palais,L-M}), so, due to automatic continuity, it extends to a
$*$-homomorphism $\underline{\Op}:C(S^*M;A)\to Q(\K\otimes A)$. Using the
Bartle--Graves selection theorem one can obtain a continuous homogeneous
lifting $\widetilde{\Op}:C(S^*M;A)\to M(\K\otimes A)$ for
$\underline{\Op}$.

Let $\gamma_0^s$ and $\gamma_{\pm 1}^s$, $s\in (0,1]$,
be smooth functions in $C_0(\mathbb R_+)$ with support in
$[2^{-1/s},2^{1/s}]$ and in $[2^{\pm 1/s-1},2^{\pm 1/s+1}]$ respectively,
such that $\sum_{i\in\mathbb Z}\gamma_i^2=1$, where
$\gamma_{\pm i}^s=\tau_{2^{\pm(i-1)}}(\gamma_{\pm 1}^s)$ for $i>1$.

Let at first $a\in C^\infty(S^*M;A)$. Define a map from $C^\infty(S^*M;A)$
to $M(\K\otimes A)$ by
 $$
\Psi_s(a)=
\sum_{i,j\in\mathbb Z}T_{1}
(\gamma_i^s\gamma_j^s\theta)\otimes a)\otimes e_{ij}
 $$
for $s\in(0,1]$
and
 $$
\Psi_0(a)=\Op(a)\otimes e_{00}.
 $$

Strict continuity of the family $\Psi_s(a)$ at any
$s\in(0,1]$ is obvious, so we have to check it at $s=0$. By construction,
$\gamma_0^s(x)=1$ for $x\in[2^{-1/s+1},2^{1/s-1}]$, hence
$\gamma_i^s$ strictly converges to zero as $s\to 0$ for
any $i\neq 0$, so for any $f\in L^2(U)\otimes A$ in local coordinates one
has
 $$
\lim_{s\to 0}\int a(x,\xi)\gamma_i^s(|\xi|)\gamma_j^s(|\xi|)\theta(|\xi|)
e^{ix\xi}\widehat{f}(\xi)\,d\xi=\left\lbrace\begin{array}{cc}
1, &{\rm if\ } i=j=0;\\0,& {\rm otherwise}\end{array}\right.
 $$
because $\widehat{f}\in L^2(\mathbb R^n)\otimes A$, hence
 \begin{equation}\label{equ1}
\lim_{s\to 0}T_1((\gamma_0^s)^2\theta\otimes a)(f)-\Op(a)(f)=0
 \end{equation}
and
 \begin{equation}\label{equ2}
\lim_{s\to 0}T_1(\gamma_i^s\gamma_j^s\theta\otimes a)(f)=0
 \end{equation}
whenever either $i$ or $j$ differs from zero. Since the set $\Psi_s(a)$
is uniformly bounded for any $a\in C^\infty(S^*M;A)$, it follows from
(\ref{equ1}) and (\ref{equ2}) that
the family of maps $\Psi_s$, $s\in[0,1]$, is strictly continuous
with respect to $s$, hence this family defines a map
 $$
\Psi:C^\infty(S^*M;A)\to M(C([0,1];\K\otimes A)),
 $$
which is obviously a $*$-homomorphism modulo the ideal $C([0,1];\K\otimes
A)$. The $*$-homomorphism $q\circ\Psi:C^\infty(S^*M;A)\to
Q(C([0,1];\K\otimes A))$ extends by continuity to a $*$-homomorphism
$\widetilde{\Psi}:C(S^*M;A)\to Q(C([0,1];\K\otimes A))$.

It remains to show that $\widetilde{\Psi}$ is the required homotopy.
One easily sees that $\widetilde{\Psi}_0=q\circ\widetilde{\Op}\oplus 0$,
so one has to check that $\widetilde{\Psi}_1=I(\widetilde{T})$ and it is
sufficient to check the latter equality on
$C^\infty(S^*M;A)$. Since, for big
enough positive $i$, $\gamma_i\gamma_j\theta=\gamma_i\gamma_j$ and since
 $$
\lim_{i\to-\infty}T_1(\gamma_i\gamma_j\otimes a)=
\lim_{i\to-\infty}T_1(\gamma_i\gamma_j\theta\otimes a)=0
 $$
for any $a\in C^\infty(S^*M;A)$ (for all $j$, because the only non-trivial
values for $j$ are $i$ and $i\pm 1$), so $q\circ \Psi_1=q\circ\Psi'$, where
 $$
\Psi'(a)=
\sum_{i,j\in\mathbb Z}T_{1}
(\gamma_i\gamma_j)\otimes a)\otimes e_{ij}.
 $$
By properties of asymptotic homomorphisms $T$ and $\widetilde{T}$ one has
 $$
\lim_{i\to-\infty}T_1(\gamma_i\gamma_j\otimes a)=0;
 $$
 $$
\lim_{i\to-\infty}\widetilde{T}_{2^i}(\tau_{2^{-i}}(\gamma_i\gamma_j)\otimes
a)=\lim_{i\to-\infty}
\widetilde{T}_{2^i}(\gamma_0\gamma_{j-i})\otimes a)=0
 $$
and
 $$
\lim_{i\to\infty}\widetilde{T}_{2^i}(\tau_{2^{-i}}(\gamma_i\gamma_j)\otimes
a)-T_1(\gamma_i\gamma_j\otimes a)
=\lim_{i\to\infty}\widetilde{T}_{2^i}(\gamma_0\gamma_{j-i}\otimes a)-
T_{2^i}(\gamma_0\gamma_{j-i}\otimes a)=0
 $$
for any $a\in C^\infty(S^*M;A)$. Therefore,
$\widetilde{\Psi}_1=q\circ\Psi'=I(\widetilde{T})$ on $C^\infty(S^*M;A)$.
\end{proof}

\section{Concluding remarks}

It may seem that the asymptotic homomorphism
$\widetilde{T}$ (\ref{asympTT}) contains
more information than the extension (\ref{extensi}) since it is defined
not only on $C_{00}(T^*M;A)$, but on the bigger $C^*$-algebra $C_0(T^*M;A)$.
In fact, the extension (\ref{extensi}) also possesses an additional
property, which is equivalent to that additional property of the
asymptotic homomorphism $\widetilde{T}$.
Namely, there is a subalgebra $C(M;A)\subset
C(S^*M;A)$ consisting of functions that are constants on the fibers and
the Busby invariant of the extension (\ref{extensi}) restricted onto
$C(M;A)$ can be lifted to $M(\K\otimes A)$. Indeed, multiplication
$\pi(a)f=af$ for $a\in C(M;A)$ and $f\in L^2(M)\otimes A$ defines such
a lifting, i.e. a $*$-homomorphism $\pi:C(M;A)\to M(\K\otimes A)$.
Using a relative version of the Bartle--Graves theorem
\cite{Loring}, one can construct a continuous section $\overline{\Op}:
C(S^*M;A)\to M(\K\otimes A)$ such that its restriction onto $C(M;A)$
coincides with $\pi$.
So we now describe how one can extend the Connes--Higson
construction to the case, when an extension of a $C^*$-algebra $D$
restricted to a $C^*$-subalgebra $C\subset D$ is liftable.
Denote the Busby invariant of such an extension by $\chi:D\to Q(B)$
and let $\overline{\chi}:D\to M(B)$ be a continuous homogeneous
lifting for $\chi$ such that $\overline{\chi}|_{C}$ is a $*$-homomorphism.
Consider the $C^*$-subalgebra
 \begin{equation}\label{alg}
C_0([0,\infty);C)\cup C_0(\mathbb R_+;D)
 \end{equation}
in $C_0([0,\infty);D)$.
The Connes--Higson construction on $SD=C_0(\mathbb R_+;D)$ can be defined
on elementary tensors of the form $f\otimes d$, $f\in C_0(\mathbb R_+)$,
$d\in D$, by the formula
 \begin{equation}\label{CH1}
CH(\chi)_t(f\otimes d)=\overline{\chi}(d)(f\circ\kappa)(u_t),
 \end{equation}
where $(u_t)_{t\in\mathbb R_+}\subset B$ is a quasicentral (with respect
to $\overline{\chi}(D)$) approximate unit, $0\leq u_t\leq 1$, and
$\kappa:(0,1]\to [0,\infty)$ is a homeomorphism (cf. \cite{CH}). In order
to extend this construction to the $C^*$-algebra (\ref{alg}) we have to
define the asymptotic homomorphism $CH(\chi)$ on
$C_0([0,\infty);C)$ compatible with (\ref{CH1}).
Let $g\otimes c\in C_0([0,\infty);C)$ be an elementary tensor,
$g\in C_0[0,\infty)$, $c\in C$. Then apply the same formula,
 $$
CH(\chi)_t(g\otimes c)=\overline{\chi}(c)(g\circ\kappa)(u_t).
 $$
In the case, when $D=C(S^*M;A)$ and $C=C(M;A)$, the $C^*$-algebra
(\ref{alg}) obviously coincides with the $C^*$-algebra $C_0(T^*M;A)$ and
the extended Connes--Higson construction gives us the asymptotic
homomorphism $\widetilde{T}$ defined on the whole $C_0(T^*M;A)$.

\vspace{2cm}
\parbox{7cm}{V. M. Manuilov\\
Dept. of Mech. and Math.,\\
Moscow State University,\\
Moscow, 119992, Russia\\
e-mail: manuilov@mech.math.msu.su
}

\end{document}